\documentclass[11pt,a4paper]{amsart}
\usepackage{amsmath, amssymb}
\usepackage{amsfonts}
\usepackage{amsthm}
\usepackage{amscd}
\usepackage{hyperref} 

\allowdisplaybreaks

\theoremstyle{plain}
\newtheorem{theorem}{{Theorem}}[subsection]

\newtheorem{proposition}[theorem]{Proposition}

\newtheorem{remark}[theorem]{Remark}
\newtheorem{definition}[theorem]{Definition}

\newtheorem{question}{Question}

\numberwithin{equation}{section}
\numberwithin{theorem}{section}
\newcount\colveccount
\newcommand*\colvec[1]{
	\global\colveccount#1
	\begin{pmatrix}
		\colvecnext
	}
	\def\colvecnext#1{
		#1
		\global\advance\colveccount-1
		\ifnum\colveccount>0
		\\
		\expandafter\colvecnext
		\else
	\end{pmatrix}
	\fi
}

\newcommand{\SL}{\mathsf{SL}}
\newcommand{\sGL}{\mathsf{GL}}
\newcommand{\fgl}{\mathfrak{gl}}

\newcommand{\bH}{\mathbb{H}}

\newcommand{\Z}{\mathbb{Z}}
\newcommand{\sad}{\mathsf{Ad}}
\newcommand{\adj}{\mathsf{adj}}
\newcommand{\tr}{\mathsf{Tr}}

\newcommand{\sHit}{\mathsf{Hit}}

\newcommand{\sT}{\mathsf{T}}

\newcommand{\sG}{\mathsf{G}}

\newcommand{\sDiag}{\mathsf{Diag}}
\newcommand{\cA}{\mathcal{A}}

\newcommand{\cX}{\mathcal{X}}
\newcommand{\cP}{\mathcal{P}}

\newcommand{\cS}{\mathcal{S}}

\newcommand{\cL}{\mathcal{L}}
\newcommand{\cY}{\mathcal{Y}}

\newcommand{\tAd}{\mathtt{Ad}}
\newcommand{\tad}{\mathtt{ad}}

\newcommand{\tB}{\mathtt{B}}

\newcommand{\tR}{\mathtt{R}}

\newcommand{\ttr}{\mathtt{tr}}

\newcommand{\tJd}{\mathtt{Jd}}
\newcommand{\tM}{\mathtt{M}}
\newcommand{\tP}{\mathtt{P}}

\newcommand{\N}{\mathbb{N}}
\newcommand{\R}{\mathbb{R}}
\newcommand{\C}{\mathbb{C}}

\newcommand{\s}{\mathsf{Stab}}

\newcommand{\sHom}{\mathsf{Hom}}

\newcommand{\defeq}{\mathrel{\mathop:}=}

\newcommand{\flow}{\mathsf{U}\Gamma}
\newcommand{\cflow}{\widetilde{\mathsf{U}\Gamma}}
\newcommand{\bdry}{\partial_\infty\Gamma}
\newcommand{\fg}{\mathfrak{g}}
\newcommand{\fa}{\mathfrak{a}}

\newcommand{\fk}{\mathfrak{k}}

\newcommand{\fp}{\mathfrak{p}}

\newcommand{\fsl}{\mathfrak{sl}}

\newcommand\irr{\operatorname{irr}}
\newcommand\mdeg{\operatorname{mdeg}}
\newcommand\mdim{\operatorname{mdim}}
\newcommand\zd{\operatorname{zd}}
\newcommand{\Spec}{\mathsf{Spec}}

\title[Affine Anosov]{Affine Anosov representations}
\author{Sourav Ghosh}
\address{Ashoka University}
\email{sourav.ghosh@ashoka.edu.in, sourav.ghosh.bagui@gmail.com}
\date{\today}

\thanks{The author acknowledge support from Ashoka University annual research grant.}
\begin{document}
	
	\begin{abstract}
		In this survey article we discuss about possible generalizations of Anosov representations in the affine setting and their consequences.
	\end{abstract}
	
	\maketitle
	\tableofcontents
	
	\newpage
	
	\section*{Introduction}
	
	The goal of this expository article is to give an overview of the study of affine Anosov representations. The notion of an affine Anosov representation arises naturally in the interface of the study of proper affine actions of hyperbolic groups and Anosov representations into semisimple Lie groups.
	
	\paragraph{\textbf{Hitchin representations}} The study of Anosov representations originates from the study of Hitchin representations. Hitchin \cite{Hit} showed that the space of representations of the fundamental group of a closed compact surface inside $\mathsf{SL}(n,\R)$ has some distinguished connected components. These components are homeomorphic to a ball of some appropriate finite dimension which depends on the genus of the surface. The components are called \emph{Hichin components} and elements of the Hitchin components are called \emph{Hitchin representations}. Moreover, for $n=2$ these connected components are classically known as Teichm\" uller spaces. Hitchin proved his result using analytic means and wondered about the geometric significance of the elements in the Hitchin components. The embeddings of convex cocompact subgroups of $\SL(2,\R)$ inside $\SL(n,\R)$ obtained via irreducible representations of $\SL(2,\R)$ inside $\SL(n,\R)$ are called \emph{Fuchsian representations}. We note that Fuchsian representations of cocompact surface groups always lie in the Hitchin components. 
	
	\paragraph{\textbf{Anosov representations}} Labourie \cite{Labourie} introduced the notion of an Anosov representation and, using results of Guichard \cite{Gui1}, gave a complete geometric description of the elements in the Hitchin components in terms of properties of the limit curve. Later, the dynamical definition of an Anosov representation given by Labourie was rephrased in Lie theoretic terms by the works of Guichard--Wienhard \cite{GW2}, Kapovich--Leeb--Porti \cite{KLPmain}, Bochi--Potrie--Sambarino \cite{BPS} and Kassel--Potrie \cite{KP} (see also Gu\' eritaud--Guichard--Kassel--Wienhard \cite{GGKW}). These works showed that the notion of an Anosov representation in some semisimple Lie group is equivalent to the image subgroup having uniform eigenvalue/singular value gaps. In some sense, Anosov subgroups are appropriate generalizations of convex cocompact subgroups in the higher rank setting (see Danciger--Gu\' eritaud--Kassel \cite{DGK4,DGK5}). Although, we must note here that the naive generalization of convex cocompactness in the higher rank setting is not interesting due to rigidity results of Mostow \cite{Mostow}, Prasad \cite{Prasad} and Margulis \cite{MargulisSuperRigid}. 
	
	\paragraph{\textbf{Crystallographic groups}}The study of proper affine actions of hyperbolic groups originates from the classification theorems of Bieberbach \cite{B1,B2} (see also the proof by Buser \cite{Buser}). Bieberbach showed that discrete subgroups $\Gamma$ of $\mathsf{O}(n,\R)\ltimes\R^n$, i.e. subgroups of rigid body motions of an $n$-dimensional Euclidean space, which act properly discontinuously, freely and cocompactly on $\R^n$ are virtually $\Z^n$. Such a group is called a \emph{crystallographic group}. Later, Auslander and Markus \cite{AM} came up with examples of discrete subgroups $\Gamma$ of $\sGL(n,\R)\ltimes\R^n$, i.e. subgroups of affine transformations of an $n$-dimensional Euclidean space, that are not isomorphic to $\Z^n$ but which act properly discontinuously, freely and cocompactly on $\R^n$. The examples constructed by Auslander and Markus are virtually polycyclic i.e. up to a finite index subgroup they are isomorphic to semidirect product of copies of $\Z$. In a failed attempt, Auslander \cite{Aus} tried to prove that any discrete subgroup $\Gamma$ of $\sGL(n,\R)\ltimes\R^n$ which act properly discontinuously, freely and cocompactly on $\R^n$ is virtually polycyclic. 
	
	\paragraph{\textbf{Margulis spacetimes}} The statement that Auslander tried to prove was later made into a conjecture, called the \emph{Auslander conjecture}. The conjecture is known to be true for lower dimensional manifolds by works of Fried--Goldman \cite{FG} and Abels--Margulis--Soifer \cite{AMS2} but it is still wide open for the general case. In hope of making the Auslander conjecture more tractable, Milnor \cite{Milnor} wondered about the possibility of removing the cocompactness assumption from the conjecture. Later, Margulis \cite{Margulis1,Margulis2} surprised everybody by producing a negative answer to Milnor's question. He showed the existence of free non-abelian subgroups of $\sGL(3,\R)\ltimes\R^3$ which act properly discontinuously and freely on $\R^3$. These examples are called \emph{Margulis spacetimes}. Subsequently, Drumm \cite{Drumm1} constructed explicit fundamental domains for a large class of the examples of Margulis. The fundamental domains constructed by Drumm are bounded by crooked planes, and it was conjectured that all examples constructed by Margulis admit fundamental domains bounded by crooked planes. After remaining open for more than a decade, the conjecture was finally resolved in the positive by Danciger--Gu\' eritaud--Kassel \cite{DGK2}.
	
	\paragraph{\textbf{Existence results}} Suppose $\Gamma\subset\sGL(3,\R)\ltimes\R^3$ is a non-abelian free group which acts properly discontinuously and freely on $\R^3$. Then by a result due to Fried--Goldman \cite{FG}, the Zariski closure of $\Gamma$ is some conjugate of $\mathsf{SO}(2,1)$. Generalising results of Margulis, Abels--Marulis--Soifer \cite{AMS} showed that there exists Zariski dense non-abelian free subgroups of $\mathsf{SO}(m-1,m)\ltimes\R^{2m-1}$ which act properly discontinuously and freely on $\R^{2m-1}$ only when $m$ is even. Later, Smilga \cite{Smilga4} constructed more such examples of proper affine actions using the arguements of Margulis more deftly (see also the work of Neza \cite{Neza} and Burelle--Neza \cite{BN}). Recently, Danciger--Gu\' eritaud--Kassel \cite{DGK3} discovered examples of Coxeter groups admitting proper affine actions. 
	\begin{question}
		What are the different kinds of discrete subgroups of the affine group $\sGL(n,\R)\ltimes\R^n$ which act properly discontinuously and freely on $\R^n$ and what are their Zariski closures?
	\end{question}
	
	\paragraph{\textbf{Non-existence results}} We note that the existence of proper affine actions also depends on the particular ways the $\Gamma$ is embedded in $\sGL(n,\R)\ltimes\R^n$. Suppose $\Gamma\subset\mathsf{SO}(2,1)\ltimes\R^3$ acts properly discontinuously and freely on $\R^3$, then Mess \cite{Mess} showed that the linear part of $\Gamma$ cannot act cocompactly on the Minkowski model of hyperbolic geometry (see also Goldman--Margulis \cite{GM}). This result was later generalized by Labourie \cite{Labourie2,Labourie3} and Danciger--Zhang \cite{DZ} to show the impossibility of proper affine actions of subgroups whose linear part is Hitchin. Using low dimensional isomorphisms of Lie groups, Mess' result can also be interpreted as stating the following:
	suppose $\Gamma\subset\SL(2,\R)\ltimes_\tAd\fsl(2,\R)$ acts properly discontinuously and freely on $\fsl(2,\R)$, then the linear part of $\Gamma$ cannot act cocompactly on the upper half plane model of hyperbolic geometry.
	Ghosh \cite{Ghosh7} presents a partial generalization of this alternative interpretation of Mess' result. 
	
	\paragraph{\textbf{Affine anosov representations}} Suppose $(g,v)\in\mathsf{SO}(2,1)\ltimes\R^3$ and $g$ is diagonalizable with distinct eigenvalues. For any such $(g,v)$, Margulis defined a number $\tM(g,v)\in\R$, which we call the \emph{Margulis invariant} of $(g,v)$, to detect proper affine actions. He showed that if $\tM(g,v)$ and $\tM(h,w)$ have opposite signs then the group generated by $(g,v)$ and $(h,w)$ cannot act properly discontinuously on $\R^3$. The notion of a Margulis invariant was later generalized by Abels--Margulis--Soifer \cite{AMS} and Smilga \cite{Smilga4}. These generalized invariants are not necessarily real valued but are vector valued quantities. Labourie \cite{Labourie2} introduced a diffused version of the Margulis invariant. Goldman, Labourie and Margulis \cite{GLM} used the positivity of these diffused Margulis invariants to give an equivalent criterion for the feasibility of proper affine actions of discrete subgroups of $\mathsf{SO}(m-1,m)\ltimes\R^{2m-1}$ whose linear parts are Fuchsian representations of convex cocompact surface groups. Taking inspiration from this work, Ghosh \cite{Ghosh2, Ghosh3, Ghosh7} and Ghosh--Treib \cite{GT} related the study of proper affine actions of word hyperbolic groups with the study of Anosov representations. In particular, they introduce the notion of an \emph{affine Anosov representation} for representations of word hyperbolic groups inside affine groups whose linear parts are semisimple Lie groups.
	
	\paragraph{\textbf{Isospectral rigidity}} We note that Margulis invariants are closely related to infinitesimal deformations of eigenvalues (see Goldman--Margulis \cite{GM}, Danciger--Zhang \cite{DZ} and Ghosh \cite{Ghosh3}). More generally, they are related to infinitesimal deformations of Jordan-Lyapunov projections (see Ghosh \cite{Ghosh7} and Sambarino \cite{Samba2}). Consequently, the marked Margulis invariant spectra of an affine representation satisfy similar kind of rigidity results that are satisfied by the marked Jordan-Lyapunov projection spectra of a linear representation (see Drumm--Goldman \cite{DG}, Charette--Drumm \cite{CD}, Kim \cite{Kim} and Ghosh \cite{Ghosh4,Ghosh5,Ghosh6}). Such rigidity results also provide nice embeddings of the corresponding representation spaces. 
	
	\paragraph{\textbf{Pressure metric}} The Teichm\" uller space has a natural complex analytic structure and a natural K\" ahler metric called the \emph{Weil--Petersson} metric \cite{Ahlfors}.  Wolpert \cite{Wolpert} showed that this metric is a multiple of the metric defined by Thurston, which was, in turn, constructed using the Hessian of the length of a random geodesic. Building on works of Bonahon \cite{Bona}, Bowen \cite{Bowen} and Bridgeman--Taylor \cite{BT}, McMullen \cite{mcmu} showed that this metric can also be interpreted in terms of the theory of thermodynamic formalism. Such metrics are called \emph{pressure metrics}. Bridgeman--Canary--Labourie--Sambarino \cite{BCLS} introduced these methods to the theory of Anosov representations and defined the pressure metric on the space of projective Anosov representations. Further properties of these metrics on the Hitchin components have been studied by Labourie--Wentworth \cite{LW} and Sambarino \cite{Samba2}. Ghosh \cite{Ghosh1,Ghosh4} used the rigidity results mentioned in the previous paragraph and applied this theory in the affine setting to construct pressure metrics on the space of affine Anosov representations into $\mathsf{SO}(2m-1,2m)\ltimes\R^{4m-1}$. Further generalization of this result to account for all affine Anosov representations is a work in progress.
	
	\subsection*{Acknowledgements}
	I would like to thank the organisers of the 2024 conference on ``Zariski dense subgroups, number theory and geometric applications", held at International Center for Theoretical Sciences, Bengaluru, for giving me the opportunity to include this article as part of the proceedings of the conference. I would also like to thank Saujanya Bharadwaj for helping me improve the exposition of this article.

	\section{A special case}
	
	Let $\fgl(2,\R)$ be the space of $2\times 2$ matrices and $\adj:\fgl(2,\R)\to\fgl(2,\R)$ be the following map:
	\[\adj\begin{bmatrix}a&b\\c&d\end{bmatrix}:=\begin{bmatrix}d&-b\\-c&a\end{bmatrix}.\]
	We consider $\langle A\mid B \rangle:= -\tr(A.\adj(B))$ and observe that it is a symmetric bilinear form of signature $(2,2)$ on $\fgl(2,\R)$. As $B$ is a $2\times 2$ matrix, the trace of $B$ is zero imply that $\adj(B)=-B$. Hence, the restriction of the symmetric bilinear form on $\fsl(2,\R)$ is a multiple of the Killing form. We consider the following action of $\mathsf{SL}(2,\R)\times\mathsf{SL}(2,\R)$ on $\fgl(2,\R)$:
	\[(g,h).A:=gAh^{-1}.\]
	We note that this action preserves the aforementioned bilinear form on $\fgl(2,\R)$. Hence, we obtain a homomorphism
	\[\varphi:\mathsf{SL}(2,\R)\times\mathsf{SL}(2,\R)\to\mathsf{O}(2,2).\]
	As $\mathsf{SL}(2,\R)$ is connected, we deduce that the image of $\varphi$ is inside $\mathsf{SO}_0(2,2)$, the connected component of $\mathsf{O}(2,2)$ containing the identity element. We observe that the image of the map $\sDiag:\mathsf{SL}(2,\R)\to\mathsf{SL}(2,\R)\times\mathsf{SL}(2,\R)$ which sends $g$ to $(g,g)$ preserves the identity matrix. Also, the signature of the Killing form on $\fsl(2,\R)$ is $(2,1)$. Hence, the image of $\sDiag(\mathsf{SL}(2,\R))$ under $\varphi$ is a copy of $\mathsf{SO}_0(2,1)$ sitting inside $\mathsf{SO}_0(2,2)$. Furthermore, we note that $\varphi$ is a double cover of $\mathsf{SO}_0(2,2)$ with kernel $\{(e,e),(-e,-e)\}$ for the identity matrix $e$.
	
	We recall that a finitely generated discrete subgroup $\Gamma$ of $\SL(2,\R)$ has a natural action on the upper half plane $\bH$ via M\" obius transformations. We call the discrete subgroup $\Gamma$ \emph{convex cocompact} if and only if $\Gamma$ acts properly discontinuously and freely on $\bH$ and the quotient manifold $\Gamma\backslash\bH$ has a compact convex core. We note that $\Gamma\backslash\bH$ has a compact convex core if and only if $\Gamma\backslash\bH$ does not have any cusps. We observe that the elements of such a group $\Gamma$ satisfy the following:
	\begin{enumerate}
		\item Any $\gamma\in\Gamma$ is diagonalizable with eigenvalues $\lambda_\gamma$ and $\lambda_\gamma^{-1}$ and $\lambda_\gamma\geq1$.
		\item There exists some $C>1$ such that for any torsion free $\gamma\in\Gamma$ the ratio of the eigenvalues satisfies $\lambda_\gamma^2\geq C$.
	\end{enumerate}
	We call the second property as a \emph{uniform eigenvalue gap property}. We consider two representations $\rho,\varrho$ of $\Gamma$ inside $\SL(2,\R)$ which are both convex cocompact. Suppose $\lambda_\rho(\gamma)$ and $\lambda_\varrho(\gamma)$ are the respective eigenvalues of the images of $\gamma\in\Gamma$ which are bounded below by $1$. We observe that the eigenvalues for the action of $(\rho(\gamma),\varrho(\gamma))$ on $\fgl(2,\R)$ are $\lambda^{\pm1}_\rho(\gamma)\lambda^{\pm1}_\varrho(\gamma)$. We observe that $$\lambda_\rho(\gamma)\lambda_\varrho(\gamma)\geq\lambda_\rho(\gamma)\lambda^{-1}_\varrho(\gamma),\lambda^{-1}_\rho(\gamma)\lambda_\varrho(\gamma)\geq\lambda^{-1}_\rho(\gamma)\lambda^{-1}_\varrho(\gamma).$$
	A priori, it is not clear which one of  $\lambda_\rho(\gamma)\lambda^{-1}_\varrho(\gamma)$ and $\lambda^{-1}_\rho(\gamma)\lambda_\varrho(\gamma)$ is bigger than the other. We observe that comparing the above two eigenvalues of $(\rho(\gamma),\varrho(\gamma))$ is equivalent to the comparison between $\lambda^2_\rho(\gamma)$ and $\lambda^2_\varrho(\gamma)$. In fact, by a result of Thurston, if both $\rho$ and $\varrho$ are cocompact (and not just convex cocompact), then neither of the numbers can be bigger than the other uniformly over all $\gamma\in\Gamma$. That is, there exists some non trivial $\gamma\in\Gamma$ such that
	$\lambda^2_\rho(\gamma)\geq\lambda^2_\varrho(\gamma)$ and some  non trivial $\eta\in\Gamma$ such that
	$\lambda^2_\rho(\eta)\leq\lambda^2_\varrho(\eta)$. 
	
	We consider a one parameter analytic family $\{\rho_t\}_{t\in(-1,1)}$ of representations $\rho_t:\Gamma\to\SL(2,\R)$ such that $\rho:=\rho_0$ is convex cocompact. Clearly, $\rho_t$ is also convex cocompact for $t$ small enough. Without loss of generality we assume that $\rho_t$ is convex cocompact for all $t\in(-1,1)$. We consider $u:\Gamma\to\fsl(2,\R)$ such that for all $\gamma\in\Gamma$,
	\[u(\gamma):=\left.\frac{d}{dt}\right|_{t=0}\rho_t(\gamma)\rho(\gamma)^{-1}.\]
	We observe that $u$ is a cocycle. We observe that the adjoint action of $\SL(2,\R)$ on its lie algebra $\fsl(2,\R)$ has a non-trivial center $\{e,-e\}$. Moreover, we also observe that the adjoint action of $\mathsf{PSL}(2,\R)$ on its lie algebra $\fsl(2,\R)$ is isomorphic to the linear action of $\mathsf{SO}_0(2,1)$ on $\R^3$. It follows that $\SL(2,\R)\ltimes_\tAd\fsl(2,\R)$ is a two fold cover of $\mathsf{SO}_0(2,1)\ltimes\R^3$. We use results of Danciger--Gu\' eritaud--Kassel \cite{DGK1} (see also Ghosh \cite{Ghosh3}) and note that
	\begin{theorem}
		Suppose $(\rho,u)(\Gamma)$ acts properly discontinuously and freely on $\fsl(2,\R)$. Then for all $t$ non-zero but close enough to zero, there exists $C_t>1$ such that the following holds:
		\begin{enumerate}
			\item either $\lambda^2_{\rho_t}(\gamma)\lambda^{-2}_{\rho}(\gamma)\geq C_t$ for all non-torsion element $\gamma\in\Gamma$,
			\item or $\lambda^2_{\rho}(\gamma)\lambda^{-2}_{\rho_t}(\gamma)\geq C_t$ for all non-torsion element $\gamma\in\Gamma$.
		\end{enumerate}
		In either case $\varphi(\rho,\rho_t)(\Gamma)\subset\mathsf{SO}_0(2,2)\subset\SL(4,\R)$ has the uniform middle eigenvalue gap property.
	\end{theorem}
	We end this section by noting that the study of Anosov representations generalize the notion of a convex cocompact representation and the study of affine Anosov representations try to capture the above theorem in more general settings.

	\section{Gromov flow space}
	
	Let $\Gamma$ be a finitely generated word hyperbolic group with generating set $S$. We denote the \emph{word length} of an element $\gamma\in\Gamma$ with respect to the generating set $S$ by $|\gamma|_S$, and the \emph{stable length} of $\gamma$ with respect to $S$ by $|\gamma|_{S,\infty}$ i.e. 
	\[    |\gamma|_S:=\min\{n\mid \gamma=s_1\dots s_n \text{ and } S\cap \{s_i,s_i^{-1}\}\neq\emptyset \text{ for all } 1\leq i\leq n\} \]
	and $|\gamma|_{S,\infty} := \lim_{n\to\infty} \frac{|\gamma^n|_S}{n}$. Let $F$ be another generating set of $\Gamma$, then $|\cdot|_F$ is bi-Lipschitz with $|\cdot|_S$ i.e. there exists some constant $C>0$ such that for all $\gamma\in\Gamma$ the following holds:
	\[\frac{1}{C}\leq \frac{|\gamma|_F}{|\gamma|_S}\leq C.\]
	It follows that the stable length corresponding to two different generating sets are also bi-Lipschitz. We note that the stable length is invariant under conjugation but the word length need not be. 
	
	We denote the \emph{Gromov boundary} of $\Gamma$ by $\bdry$. There is a natural metric on the space $\bdry$ called the \emph{visual metric}. This metric is unique up to H\" older equivalence (for more details, see Theorem 2.18 of Benakli--Kapovich \cite{BK}). 
	
	We note that the fundamental group of a surface with negative Euler characteristic is a hyperbolic group. Henceforth, we will call such a group a \emph{surface group}. When $\Gamma$ is a surface group, the Gromov boundary is homeomorphic to the limit set of the action of the surface group on the upper half plane $\bH$. 
	
	The group $\Gamma$ has a natural action on $\bdry$ and for each $\gamma\in\Gamma$ of infinite order, there exist two distinct points on $\bdry$ that are fixed by the action of $\gamma$. One of these points is an attracting fixed point of $\gamma$ and the other one is a repelling fixed point of $\gamma$. We denote the attracting fixed point of $\gamma$ by $\gamma^+$ and the repelling fixed point of $\gamma$ by $\gamma^-$. We note that the attracting fixed point of $\gamma$ is the repelling fixed point of $\gamma^{-1}$ and vice-versa. The collection of attracting fixed points of all infinite order elements $\gamma\in\Gamma$ form a dense subset of $\bdry$. For a more precise statement see Propositions 4.2 and 4.3 of Benakli--Kapovich \cite{BK}. We denote the subset of $\bdry^2$ obtained by removing the diagonal by $\bdry^{(2)}$ and consider the following set:
	\[\cflow:=\bdry^{(2)}\times\R.\]
	The visual metric on $\bdry$ naturally induces a metric on $\cflow$, called the product metric. Let $\phi_t:\cflow\to\cflow$ be the map such that $\phi_t(x,y,s)=(x,y,s+t)$ for all $t\in\R$. By a celebrated result of Gromov \cite{Gromov}, there exists a proper action of $\Gamma$ on $\cflow$ which extends the diagonal action of $\Gamma$ on $\bdry^{(2)}$ and which commutes with the maps $\{\phi_t\}_{t\in\R}$. Moreover, there exists a metric $d$ on $\cflow$ which is bi-Lipschitz to the product metric on $\cflow$ and for which the following holds:
	\begin{enumerate}
		\item $\Gamma$ is a subgroup of the isometry group of $(\cflow,d)$,
		\item $\phi_t$ are Lipschitz homeomorphisms,
		\item orbits of $\phi_t$ are quasi-isometrically embedded in $\cflow$.
	\end{enumerate}
	The \emph{Gromov flow space} of $\Gamma$, denoted by $\flow$, is the quotient space of $\cflow$ under the action of $\Gamma$. We abuse notation and denote the maps on $\flow$ induced from the maps $\phi_t$ by $\phi_t$. 
	
	We note that when $\Gamma$ is a surface group, $\bdry^{(2)}$ identifies with the space of all bi-recurrent geodesics of $\Gamma\backslash\bH$ and the Gromov flow space identifies with the following:
	\[\{(x,l)\mid l \text{ is a bi-recurrent geodesic of }\Gamma\backslash\bH \text{ and } x\in l\}.\]
	Moreover, when the hyperbolic surface $\Gamma\backslash\bH$ has no boundary, the space $\cflow$ identifies with $\mathsf{PSL}(2,\R)$ and the Gromov flow space identifies with $\Gamma\backslash\mathsf{PSL}(2,\R)$ which in turn identifies with the unit tangent bundle of the surface $\Gamma\backslash\bH$.
	
	The space $\flow$ is a connected proper metric space (for more details see Champetier \cite{Champetier} and Mineyev \cite{Mineyev}). Suppose $\gamma\in\Gamma$ is an element of infinite order. Then the following limit is independent of $p\in\cflow$ and is called the \emph{translation length} of $\gamma$:
	\[\ell(\gamma):=\lim_{n\to\infty}\frac{d(\gamma^np,p)}{n}.\]
	The translation length of $\gamma$ is also the infimum of the distance between any point $p\in\flow$ and $\gamma p$. Moreover, the infimum is attained when $p$ is of the form $(\gamma^-,\gamma^+,t)$ for some $t\in\R$. 
	
	We note that when $\Gamma$ is a surface group, $\ell(\gamma)$ is equal to the length of the closed geodesic corresponding to $\gamma$.
	
	\section{Anosov representations}
	
	Let $G$ be a real semisimple Lie group and let $P^\pm$ be a pair of opposite parabolic subgroups of $G$. Hence, $G/P^\pm$ is necessarily compact. Let $\Gamma$ be a finitely generated hyperbolic group and $\rho:\Gamma\to G$ be an injective homomorphism. We note that $G$ acts naturally on $G/P^\pm$ as follows: for any $g\in G$ and $[h]:=hP^\pm\in G/P^\pm$ we have $g[h]:=[gh]$. We observe that under the diagonal action of $G$ on $(G/P^+)\times(G/P^-)$ there exists a unique open orbit. We denote it by $\cX$ and note that $\cX$ naturally identifies with the quotient $G/(P^+\cap P^-)$. Hence, for any $x\in\cX$, $x$ identifies with $(gP^+,gP^-)$ for some $g\in G$. We denote $gP^\pm$ respectively by $x^\pm$.
	
	We consider the bundle $\cflow\times \cX$ over $\cflow$. The map $\phi_t:\cflow\to\cflow$ induces a map $\phi_t: \cflow\times \cX \to \cflow\times \cX$ which takes $(p,x)$ to $(\phi_t p,x)$.  The group $\Gamma$ acts on $\cflow\times \cX$ by taking $(p,x)$ to $(\gamma p, \rho(\gamma)x)$ for all $p\in\cflow$ and $x\in\cX$. We take the quotient under this action of $\Gamma$ and obtain a bundle $\cX_\rho$ over $\flow$. The map $\phi_t$ again induces a map $\phi_t:\cX_\rho\to\cX_\rho$. We note that by construction $\cX$ comes equipped with a pair of distributions $X^\pm$ such that $(X^\pm)_{x}=\sT_{x^\pm}(G/P^\pm)$. As the action of $G$ on $X^\pm$ keeps them invariant, we observe that $X^\pm$ can be naturally interpreted as vector bundles over $\cX_\rho$.
	
	\begin{definition}  \label{Def:Anosov} 
		A representation $\rho:\Gamma\to G$ is called $(P^+,P^-)$-\emph{Anosov} if there exists a section $\sigma: \flow \to \cX_\rho$ such that
		\begin{itemize}
			\item $\sigma$ is locally constant along the flow lines of the flow $\phi_t$, with respect to the locally flat structure on $\cX_\rho$,
			\item the flow $\phi_t$ is contracting on the bundle $\sigma^*X^+$ and dilating on the bundle $\sigma^*X^-$.
		\end{itemize}
	\end{definition}
	
	We note that existence of a $\sigma: \flow \to \cX_\rho$ which is locally constant along the flow lines of the flow $\phi_t$ is equivalent to the existence of a map $\sigma:\cflow\to \cX$ such that $\sigma(\gamma p)=\rho(\gamma)\sigma(p)$ and $\sigma(\phi_t p)=\sigma(p)$ for all $\gamma\in\Gamma$, $p\in\cflow$ and $t\in\R$. Moreover, to say that the flow $\phi_t$ is contracting on the bundle $\sigma^*X^+$ and dilating on the bundle $\sigma^*X^-$ is equivalent to saying that there exists a pair of positive constants $c,C$ and for all $p\in\cflow$ there exists a collection of Euclidean norms $\|\cdot\|_{p}$ on $\sT_{\sigma(p)}\cX$ such that for all $\gamma\in\Gamma$ and $v\in\sT_{\sigma(p)}\cX$ we have $\|\rho(\gamma)v\|_{\gamma p}=\|v\|_p$ and for all $t>0$ and $v_\pm\in\sT_{\sigma(p)^\pm}(G/P^\pm)=\sT_{\sigma(\phi_t p)^\pm}(G/P^\pm)$ we have $\|v_\pm\|_{\phi_{\pm t} p}\leq Ce^{-ct}\|v_\pm\|_p$. 
	
	Suppose $G$ is $\SL(n,\R)$. We consider the following two subspaces of $\R^n$: $\R^k\times\{0\}^{n-k}$ and $\{0\}^k\times\R^{n-k}$. We observe that they are transverse to each other. Suppose $P_k^+$ is the stabilizer of the aforementioned $k$-dimensional subspace and $P_k^-$ is the stabilizer of the aforementioned $(n-k)$-dimensional subspace. We note that $P_k^\pm$ form a pair of opposite parabolic subgroups of $\SL(n,\R)$. 
	
	We recall that for any element $g\in\SL(n,\R)$, $gg^t$ is a symmetric matrix. Hence, $\sqrt{gg^t}$ is well defined. The eigenvalues of $\sqrt{gg^t}$ is called the \emph{singular values} of $g$. The singular values are necessarily positive. Let $\{\kappa_{i}(g)\}_{i=1}^n$ be the logarithm of the singular values of $g$ such that $\kappa_1(g)\geq\dots\geq\kappa_n(g)$. We denote the tuple $(\kappa_1(g),\dots,\kappa_n(g))$ by $\kappa(g)$. Moreover, suppose $\{\lambda_{i}(g)\}_{i=1}^n$ be the logarithm of the modulus of the eigenvalues of $g$ such that $\lambda_1(g)\geq\dots\geq\lambda_n(g)$. We denote the tuple $(\lambda_1(g),\dots,\lambda_n(g))$ by $\lambda(g)$. We note that for any $g\in G$ the following identity holds:
	\[\lim_{n\to\infty}\frac{\kappa(g^n)}{n}=\lambda(g).\]
	
	Let $\rho:\Gamma\to G$ be an injective homomorphism. We say that $\rho$ has a uniform gap at its $i$-th singular value if there exist $c,k>0$ such that for all $\gamma\in\Gamma$ the following holds:
	\[\kappa_i(\rho(\gamma))-\kappa_{i+1}(\rho(\gamma))\geq c|\gamma|_S-k.\]
	Moreover, we say that $\rho$ has a uniform gap at its $i$-th eigenvalue value if there exists $c,k>0$ such that for all $\gamma\in\Gamma$, the following holds:
	\[\lambda_i(\rho(\gamma))-\lambda_{i+1}(\rho(\gamma))\geq c|\gamma|_{S,\infty}-k.\]
	We note that having a uniform gap at its $i$-th singular value or having a uniform gap at its $i$-th eigenvalue is independent of the choice of the generating set. 
	\begin{theorem}[\cite{KLPmain}, \cite{BPS}]
		Suppose $\rho:\Gamma\to \SL(n,\R)$ is an injective homomorphism and $\Gamma$ is word hyperbolic. Then $\rho$ is Anosov with respect to the parabolic subgroups $P_i^\pm$ if and only if $\rho$ has a uniform gap at its $i$-th singular value.
	\end{theorem}
	
	\begin{theorem}[\cite{KP}]
		Suppose $\rho:\Gamma\to \SL(n,\R)$ is an injective homomorphism and $\Gamma$ is word hyperbolic. Then $\rho$ has a uniform gap at its $i$-th singular value if and only if $\rho$ has a uniform gap at its $i$-th eigenvalue.
	\end{theorem}

	\section{Proper affine actions of free groups}
	
	Suppose $G$ is a semisimple Lie group without compact factors. We denote its Lie algebra by $\fg$. Suppose $C_g:G\to G$ be the conjugation action which sends $h\in G$ to $ghg^{-1}$. We denote the identity element of $G$ by $e$ and, the differential of $C_g$ at $e$ by $\tAd(g)$ i.e.  $\tAd(g):=d_eC_g$. Hence, $\tAd$ is a map from $G$ to $\SL(\fg)$. We denote $d_e\tAd$ by $\tad$ and note that $\tad:\fg\to\fsl(\fg)$. Suppose $X,Y\in\fg$. We define $\tB(X,Y):=\tr(\tad(X)\tad(Y))$. We note that $\tB$ is a non-degenerate symmetric bilinear form on $\fg$ and is called the \emph{Killing form}. Suppose $\theta:\fg\to\fg$ be such that $\theta^2$ is trivial and $-\tB(\cdot,\theta \cdot)$ is positive definite, then $\theta$ is called a \emph{Cartan involution} of $\fg$. We note that every real semisimple Lie algebra admits a Cartan involution (for more details see Corollary 6.18 of Knapp \cite{Knapp}). We denote the eigenspace of $\theta$ of eigenvalue $1$ (resp. $-1$) by $\fk$ (resp. $\fp$). We note that $\fg=\fk\oplus\fp$. We denote the maximal abelian subspace of $\fp$ by $\fa$ and denote its dual space by $\fa^*$ i.e. $\fa^*$ is the space of all linear forms on $\fa$. We denote the connected subgroup of $G$ whose Lie algebra is $\fa$ by $A$ and the centralizer of $A$ in $G$ by $L$ i.e. $L=\{g\in G\mid C_g(a)=a \text{ for all } a\in A\}$. We consider
	\[W:=\{g\in G\mid C_g(A)=A\}/L.\]
	The group $W$ is called the \emph{restricted Weyl group}. Suppose $\alpha\in\fa^*$ and
	\[\fg^\alpha:=\{X\in\fg\mid \tad(H)(X)=\alpha(H)X \text{ for all } H\in\fa\}.\]
	If $\fg^\alpha\neq0$ and also $\alpha\neq0$, then $\alpha$ is called a \emph{restricted root}. We denote the set of all restricted roots by $\Sigma$ and note that $\Sigma$ is finite and
	\[\fg=\fg^0\bigoplus_{\alpha\in\Sigma}\fg^\alpha.\]
	Suppose $\alpha$ is a restricted root. We observe that $(\tAd_g)^*(\alpha):=\alpha\circ\tAd_g$ is also a restricted root for all $g\in G$ normalizing $A$. Hence, the restricted Weyl group has a natural action on the collection of restricted roots. The hyperplanes $\{\ker(\alpha)\}_{\alpha\in\Sigma}$ divide $\fa$ into finitely many connected components. We choose a connected component of $\fa\setminus\bigcup_{\alpha\in\Sigma}\ker(\alpha)$ and denote it by $\fa^{++}$ and its closure by $\fa^+$. We define
	\[\Sigma^+:=\{\alpha\in\Sigma\mid\alpha(H)\geq0 \text{ for all } H\in\fa^+\}.\]
	Any element of $\Sigma^+$ is called a positive root. We note that if $\alpha$ is a restricted root, then $-\alpha$ is also a restricted root. Suppose $-\Sigma^+$ be the collection of all $\alpha\in\Sigma$ such that $-\alpha$ is a positive root. We note that $\Sigma=\Sigma^+\cup -\Sigma^+$ and $\Sigma^+\cap -\Sigma^+=\emptyset$. Moreover, we note that there exists a unique element $w_0\in W$ such that $w_0(\Sigma^+)=-\Sigma^+$. The element $w_0$ is called the \emph{longest element} of the restricted Weyl group (for reasons which we will not go into in this article).  
	
	\begin{remark}
		A word of caution: observe that although $w_0(\Sigma^+)=-\Sigma^+$, it does not necessarily imply that $w_0(\alpha)=-\alpha$ for $\alpha\in\Sigma^+$. In fact, it is not true in general.
	\end{remark}
	
	\begin{theorem}[\cite{Smilga4}]
		Suppose $G$ is a real semisimple Lie group which is connected. Suppose $V$ is a finite dimensional real vector space, $\tR:G\to\SL(V)$ is an injective homomorphism, and there exists a $v\in V$ such that the following conditions hold:
		\begin{enumerate}
			\item $lv=v$ for all $l\in L$,
			\item $\tilde{w}_0v\neq v$ for some (any) $\tilde{w}_0\in G$ with $\tilde{w}_0L=w_0\in W$.
		\end{enumerate}
		Then there exists a free non-abelian subgroup $\Gamma$ of $G\ltimes_\tR V$ whose linear part is Zariski dense in $G$ and which acts properly discontinuously on the affine space corresponding to $V$.
	\end{theorem}
	
	We note that in particular this result holds when $G=\mathsf{SO}(2n,2n-1)$ and $\tR$ is the natural inclusion of $\mathsf{SO}(2n,2n-1)$ inside $\SL(4n-1,\R)$ which is a result by Abels--Margulis--Soifer \cite{AMS}. These results are generalizations of a celebrated result due to Margulis \cite{Margulis1, Margulis2} which shows the existence of non-abelian free subgroups of $\SL(3,\R)\ltimes\R^3$ that act properly discontinuously and freely on $\R^3$.
	
	We end this section with the following negative result:
	\begin{theorem}[\cite{AMS}]
		Suppose $G=\mathsf{SO}(2n+1,2n)$ and $\tR$ is the natural inclusion of $\mathsf{SO}(2n+1,2n)$ inside $\SL(4n+1,\R)$. Then there does not exist any subgroup $\Gamma$ of $G\ltimes_\tR V$ whose linear part is Zariski dense inside $G$ and which acts properly discontinuously on $V$.
	\end{theorem}

	\section{Equivalent criterion for hyperbolic groups}
	
	Suppose $G$ is a real split Lie group, $V$ is a finite dimensional real vector space of dimension at least two and $\tR:G\to\SL(V)$ is a faithful irreducible algebraic representation. We consider its differential $d_e\tR:\fg\to\fsl(V)$, and for $\omega\in\fa^*$, we define
	\[V^\omega:=\{v\in V\mid d_e\tR(H)(v)=\omega(H)v \text{ for all } H\in\fa \}.\]
	Any $\omega\in\fa^*$ for which $V^\omega\neq0$ is called a \emph{restricted weight} of the representation $\tR$. We denote the set of all restricted weights of $\tR$ by $\Omega$. We note that $\Omega$ is a finite set and $V=\bigoplus_{\omega\in\Omega}V^\omega$. We consider the following for any $X\in\fa$:
	\[\Omega^\pm(X):=\{\omega\in\Omega\mid\pm\omega(X)>0\} \text{ and } \Omega^{0}(X):=\{\omega\in\Omega\mid\omega(X)=0\}.\]
	We say two elements $X,Y$ of $\fa$ are of the \emph{same type}, denoted by $X\sim Y$, if $\Omega^+(X)=\Omega^+(Y)$ and $\Omega^-(X)=\Omega^-(Y)$. Any element $X\in\fa$ is called
	\begin{enumerate}
		\item \emph{generic}, if $\Omega^0(X)\subset\{0\}$,
		\item \emph{symmetric}, if $w_0(X)=-X$,
		\item \emph{extreme}, if $\{w\in W\mid wX\sim X\}=\{w\in W\mid wX= X\}$.
	\end{enumerate}
	A representation $\tR$ is called \emph{non-swinging} if there exists $X\in\fa$ such that $X$ is generic, symmetric and extreme. Henceforth, we only consider non-swinging representations. We fix a generic, symmetric and extreme element of $\fa^+$ and denote it by $X_\tR$ (for more details see Smilga \cite{Smilga3}). We fix the following subspaces of $V$:
	\[V^\pm:=\bigoplus_{\omega\in\Omega^\pm(X_\tR)}V^\omega \text{ and } V^0:=\bigoplus_{\omega\in\Omega^0(X_\tR)}V^\omega,\]
	and observe that $V=V^+\oplus V^0 \oplus V^-$. Moreover, we denote $V^{\pm,0}:=V^\pm\oplus V^0$.
	
	We note that certain aspects of the action of $G$ on $V$ can be read off from its adjoint action on $\fg$. Indeed, we consider the following subspaces of $\fg$:
	\[\fp^{\pm,0}:=\fg^0\bigoplus_{\{\alpha\in\Sigma\mid \pm\alpha(X_\tR)\geq0\}}\fg^\alpha.\]
	Let $P_\tR^\pm:=\{g\in G\mid C_g(\fp^{\pm,0})=\fp^{\pm,0}\}$ be the normalizer of $\fp^{\pm,0}$ inside $G$ and $\s_G(W):=\{g\in G\mid \tR(g)W\subset W\}$ be the stabilizer of any subspace $W$ of $V$. We note that (for more details see Smilga \cite{Smilga3})
	\[P_\tR^\pm=\s_G(V^{\pm,0})=\s_G(V^{\pm}).\]
	In this section, we will only consider injective homomorphisms $\rho:\Gamma\to G$ such that $\rho$ is Anosov in $G$ with respect to these parabolic subgroups $P_\tR^\pm$.
	
	We recall that $\fk$ is the eigenspace of eigenvalue $1$ for the Cartan involution $\theta:\fg\to\fg$. We denote the connected subgroup whose Lie algebra is $\fk$ by $K$ and note that it is a maximal compact subgroup of $G$. We denote $A^+:=\exp(\fa^+)$. Moreover, we denote the connected subgroup of $G$ whose Lie algebra is $\bigoplus_{\alpha\in\Sigma^+}\fg^\alpha$ by $N$. We note that any element $g\in G$ can be written as a product of three unique elements $g_e,g_h,g_u\in G$ such that $g_e$ is conjugate to some element in $K$, $g_h$ is conjugate to some element in $A$, $g_u$ is conjugate to some element in $N$ and $g_e, g_h, g_u$ commute with each other. The subscripts are chosen keeping in mind that elements conjugate to elements in $K$ are called \emph{elliptic},  elements conjugate to elements in $A$ are called \emph{hyperbolic} and elements conjugate to elements in $N$ are called \emph{unipotent}. Suppose $g\in G$ with hyperbolic part $g_h$. Then the \emph{Jordan projection} of $g$, denoted by $\tJd_g$, is the unique element in $\fa^+$ such that $g_h$ is conjugate to $\exp(\tJd_g)$. We note that when $G=\SL(n,\R)$, the Jordan projection of an element $g\in G$ is the same as the ordered tuple $\lambda(g)$ of the logarithm of the modulus of the eigenvalues of $g$. 
	
	Now we introduce certain invariants which were introduced originally by Margulis \cite{Margulis1,Margulis2} to have more control on detecting proper actions.
	
	\begin{definition}
		Suppose $(g,v)\in G\ltimes_\tR V$ is such that $\tJd_g$ is of the same type as $X_\tR$ and $h\in G$ is such that the hyperbolic part $g_h=h\exp(\tJd_g)h^{-1}$. Also, suppose $\pi_0: V\to V^0$ is the projection map with respect to the splitting $V=V^+\oplus V^0\oplus V^-$. Then the \emph{Margulis invariant} of $(g,v)$ is defined as:
		\[\tM(g,v):=\pi_0(h^{-1}v).\]
	\end{definition}
	
	We end this section with a recent result showing the relation between Margulis invariants and proper affine actions. In order to do that, we define the notion of a diverging sequence. Suppose $\{\gamma_n\}_{n\in\N}$ is a sequence of elements in $\Gamma$ such that the word length of $\gamma_n$ goes to infinity as $n$ goes to infinity and also $\lim_{n\to\infty}\gamma_n^+=a\neq b=\lim_{n\to\infty}\gamma_n^-$. Then we call such a sequence a \emph{diverging sequence}. We note that the translation length of a diverging sequence diverges. We also note that sequences of the kind $\{\gamma^n\eta\gamma^{-n}\}_{n\in\N}$, are not diverging sequences according our notion, although their word length diverges.
	
	\begin{theorem}[\cite{Ghosh7}]
		Suppose $\Gamma$ is word hyperbolic, $\tR$ is non-swinging and $\rho:\Gamma\to G$ is an injective homomorphism which is Anosov in $G$ with respect to the parabolic subgroups $P^\pm_\tR$. Suppose $u:\Gamma\to V$ is a cocycle with respect to $\tR\circ\rho$, i.e. $(\rho,u)(\Gamma)$ is a subgroup of $G\ltimes_\tR V$.
		Then the action of $(\rho,u)(\Gamma)$ on $V$ is not proper if and only if there exists a diverging sequence $\{\gamma_n\}_{n\in\N}$ of elements of $\Gamma$ such that ${\tM(\rho(\gamma_n),u(\gamma_n))}$ stays bounded. 
	\end{theorem}
	
	\section{Representations of surface groups}
	
	In this section we show that the existence of proper actions depends not only on the particular linear representation of the semisimple Lie group but also on the representation of the hyperbolic group inside the semisimple Lie group.
	
	Suppose $\Gamma$ is the fundamental group of a compact orientable surface of negative Euler characteristic without boundary. Hence, there exists natural injective homomorphisms of $\Gamma$ inside $\SL(2,\R)$ which relate to complex/hyperbolic structures on the surface. Subgroups of $\SL(2,\R)$ obtained from such homomorphisms are called \emph{Fuchsian groups}. We recall that there exists a unique (up to conjugacy) irreducible representation of $\fsl(2,\R)$ on a finite dimensional vector space. Let $\iota:\SL(2,\R)\to\SL(n,\R)$ be the corresponding homomorphism. We say that a representation $\rho:\Gamma\to\SL(n,\R)$ is \emph{Fuchsian} if it is a composition of a representation into $\SL(2,\R)$ and the irreducible representation $\iota$. We denote the space of all homomorphisms of $\Gamma$ into $\SL(n,\R)$ by $\sHom(\Gamma,\SL(n,\R))$. We observe that for any $\rho$ in $\sHom(\Gamma,\SL(n,\R))$ and $g$ in $\SL(n,\R)$, $C_g\circ\rho$ is also a representation. We consider the quotient of the space of all representations under this action of $\SL(n,\R)$ on $\sHom(\Gamma,\SL(n,\R))$ and denote the resulting space by $\sHom(\Gamma,\SL(n,\R))/\SL(n,\R)$. We call the connected components of the representation space $\sHom(\Gamma,\SL(n,\R))/\SL(n,\R)$, which contain Fuchsian representations as \emph{Hitchin components}. Due to a celebrated result of Hitchin \cite{Hit}, these connected components are topologically trivial i.e. homeomorphic to a finite dimensional ball. We call representations lying in the Hitchin components as \emph{Hitchin representations} and denote the space of all Hitchin representations by $\sHit_n$. Due to an important result of Labourie \cite{Labourie} we know that any Hitchin representation is Anosov with respect to a pair of opposite minimal parabolic subgroups of $\SL(n,\R)$. We note that the property of being an Anosov representation is not a unique feature of representations only in the Hitchin component. Indeed, there are representations in $\SL(3,\R)$, due to Barbot \cite{Barbot}, which are Anosov with respect to a pair of opposite minimal parabolic subgroups but which do not lie in the Hitchin components. 
	
	\begin{theorem}[\cite{DZ,Labourie3}]
		Suppose $\rho:\Gamma\to\SL(n,\R)$ is a Hitchin representation and $u:\Gamma\to\R^n$ is a cocycle with respect to $\rho$. Then $(\rho,u)(\Gamma)\subset\SL(n,\R)\ltimes\R^n$ does not act properly discontinuously on $\R^n$.
	\end{theorem}
	
	The above theorem is a generalization of important results due to Mess \cite{Mess}, Goldman--Margulis \cite{GM} and Labourie \cite{Labourie2}. Moreover, we note that in the above theorem both the following facts are crucial:
	\begin{enumerate}
		\item $\Gamma$ is the fundamental group of a compact orientable surface of negative Euler characteristic without boundary,
		\item  $\rho:\Gamma\to\SL(n,\R)$ is a Hitchin representation.
	\end{enumerate}
	In fact, if one does not assume that the representation is Hitchin, then as a particular case of a result by Danciger--Gu\'eritaud--Kassel \cite{DGK3}, we know that the fundamental group of a compact orientable surface of negative Euler characteristic without boundary does admit proper affine actions. 
	
	Suppose $\sHit_n$ is the space of all Hitchin representations. We consider $\SL(n,\R)\ltimes\fsl(n,\R)$ where $\SL(n,\R)$ acts on $\fsl(n,\R)$ via the adjoint representation. We note that Mess' result \cite{Mess} can also be interpreted as follows: if $(\rho,u):\Gamma\to\SL(2,\R)\ltimes\fsl(2,\R)$ is an injective homomorphism with $\rho\in\sHit_2$, then $(\rho,u)(\Gamma)$ does not act properly on $\fsl(2,\R)$. This result naturally leads us to the following question:
	\begin{question}\label{qn.2}
		Suppose $(\rho,u):\Gamma\to\SL(n,\R)\ltimes\fsl(n,\R)$ is an injective homomorphism with $\rho\in\sHit_n$. Is it true that for all cocycles $u$, $(\rho,u)(\Gamma)$ does not act properly on $\fsl(n,\R)$? If the answer is negative, then can we classify all such cocycles $u$ for which $(\rho,u)(\Gamma)$ act properly on $\fsl(n,\R)$?
	\end{question}
	
	We end this section with a result, which provides a partial answer to Question \ref{qn.2}. We note that $\sHit_n$ is an analytic manifold (for more details see Bridgeman--Canary--Labourie--Sambarino \cite{BCLS}). Moreover, all representations contained in $\sHit_2$ are by definition Fuchsian. Suppose $\rho\in\sHit_n$ is a Fuchsian representation. We note that the tangent space to $\sHit_n$ at $\rho$ identifies with the direct sum of $\sT_k$, the spaces of $k$-differentials, for $2\leq k\leq n$. In local coordinates a $k$-differential is the symmetric product of $k$-many  differential $1$-forms (for more details see Labourie--Wentworth \cite{LW} and Sambarino \cite{Samba2}). Suppose $\sT_{odd}$ denote the direct sum of $\sT_k$ for all odd $k$ and $\sT_{even}$ denote the direct sum of $\sT_k$ for all even $k$. Clearly, $\sT_\rho\sHit_n=\sT_{odd}\oplus\sT_{even}$. 
	\begin{theorem}
		Suppose $(\rho,u):\Gamma\to\SL(n,\R)\ltimes\fsl(n,\R)$ is an injective homomorphism such that $\rho$ is Fuchsian and $u\in\sT_{odd}\oplus\sT_{2m}$ for some $1\leq m\leq \lfloor\frac{n}{2}\rfloor$. Then $(\rho,u)(\Gamma)$ does not act properly on $\fsl(n,\R)$.
	\end{theorem}
	
	\section{Affine versions of Anosov representations}
	
	We recall that the notion of an Anosov representation is equivalent to having uniform singular value or eigenvalue gaps. We take inspiration from this definition and also closely notice the equivalent criterion for proper affine actions in terms of Margulis invariants to define the notion of an affine Anosov representation. 
	
	Suppose $\rho:\Gamma\to G$ is an injective homomorphism which is Anosov in $G$ with respect to the parabolic subgroups $P^\pm_\tR$ and $u:\Gamma\to V$ is a cocycle i.e. $(\rho,u)(\Gamma)$ is a subgroup of $G\ltimes_\tR V$. We define the normalized Margulis invariant spectra for such a representation as follows:
	\[\tM\text{-}\Spec(\rho,u):=\overline{\left\{\frac{\tM(\rho(\gamma),u(\gamma)}{\ell(\gamma)}\mid \gamma\in\Gamma\right\}}.\]
	We note that $\tM\text{-}\Spec(\rho,u)$ is a convex set (for more details see Ghosh \cite{Ghosh7} or Sambarino \cite{Samba2}). We also observe that if the action of $(\rho,u)(\Gamma)$ on $V$ is not proper, then there exists a sequence $\{\gamma_n\}_{n\in\N}$ of elements of $\Gamma$ such that ${\tM(\rho(\gamma_n),u(\gamma_n))}$ stays bounded and $\{\ell(\gamma_n)\}_{n\in\N}$ diverges. Hence, it follows that $0\in\tM\text{-}\Spec(\rho,u)$.
	
	We consider the stabilizer of $V^{\pm,0}$ inside $G\ltimes_\tR V$ and respectively denote it by $Q_\tR^\pm$. We note that $G\ltimes_\tR V$ has a natural action on the coset space $(G\ltimes_\tR V)/Q_\tR^\pm$. We observe that under the diagonal action of $G\ltimes_\tR V$ on  $(G\ltimes_\tR V)/Q_\tR^+ \times(G\ltimes_\tR V)/Q_\tR^-$, there exists a unique open orbit. We denote it by $\cY$ and note that it identifies with the quotient $(G\ltimes_\tR V)/(Q_\tR^+\cap Q_\tR^-)$. Hence, for any $y\in\cY$, $y$ identifies with $((g,v)Q^+,(g,v)Q^-)$ for some $(g,v)\in G\ltimes_\tR V$. We denote $(g,v)Q^\pm$ respectively by $y^\pm$. 
	
	Similarly as before, we consider the bundle $\cflow\times\cY$ over $\cflow$ and the quotient of this bundle under the action of $\Gamma$ to obtain the bundle $\cY_{\rho,u}$ over $\flow$. As before, we obtain a natural map $\phi_t:\cY_{\rho,u}\to\cY_{\rho,u}$ extending the map $\phi_t$ on $\cflow$. We note that by construction $\cY$ comes equipped with a pair of distributions $Y^\pm$ such that $(Y^\pm)_{x}=\sT_{y^\pm}((G\ltimes_\tR V)/Q_\tR^\pm)$. As the action of $G\ltimes_\tR V$ on $Y^\pm$ keeps them invariant, we observe that $Y^\pm$ can be naturally interpreted as vector bundles over $\cY_{\rho,u}$.
	\begin{definition} 
		An injective homomorphism $(\rho,u):\Gamma\to G\ltimes_\tR V$ is called \emph{partially} $(Q_\tR^+,Q_\tR^-)$-\emph{affine Anosov} if there exists a section $\sigma: \flow \to \cY_{\rho,u}$ such that
		\begin{itemize}
			\item $\sigma$ is locally constant along the flow lines of the flow $\phi_t$, with respect to the locally flat structure on $\cY_{\rho,u}$,
			\item the flow $\phi_t$ is contracting on the bundle $\sigma^*Y^+$ and dilating on the bundle $\sigma^*Y^-$.
		\end{itemize}
	\end{definition}
	
	As before, we note that the first condition is equivalent to the existence of a $(\rho,u)(\Gamma)$-equivariant map $\sigma:\cflow\to \cY$ which is also invariant under the flow. Moreover, the second condition is equivalent to the existence of a pair of positive constants $c,C$ and for all $p\in\cflow$ a collection of Euclidean norms $\|\cdot\|_{p}$ on $\sT_{\sigma(p)}\cY$ such that for all $\gamma\in\Gamma$ and $v\in\sT_{\sigma(p)}\cY$ we have $\|\rho(\gamma)v\|_{\gamma p}=\|v\|_p$ and for all $t>0$ and $v_\pm\in\sT_{\sigma(p)^\pm}((G\ltimes_\tR V)/Q_\tR^\pm)=\sT_{\sigma(\phi_t p)^\pm}((G\ltimes_\tR V)/Q_\tR^\pm)$ we have $\|v_\pm\|_{\phi_{\pm t} p}\leq Ce^{-ct}\|v_\pm\|_p$. 
	
	We note that an affine representation is called partially affine Anosov even when it satisfies seemingly all the properties listed in the definition of an Anosov representation. The reason becomes apparent from the next proposition which relates the Anosov property of $\rho$ with the affine Anosov property of $(\rho,u)$.
	
	\begin{proposition}
		Suppose $\Gamma$ is word hyperbolic, $\tR$ is non-swinging and $(\rho,u):\Gamma\to G\ltimes_\tR V$ is an injective homomorphism. Then $(\rho,u)$ is partially affine Anosov with respect to $Q_\tR^\pm$ if and only if $\rho$ is Anosov with respect to $P_\tR^\pm$.
	\end{proposition}
	
	We note that uniform eigenvalue (resp. singular value) gap is a central feature of an Anosov representation. In the affine case, the appropriate analogue is a condition on the spectra of Margulis invariants. Appropriate analogues of Anosov representations in the affine case first appeared in \cite{Ghosh1,Ghosh2}. The following definition is a natural extension of that.
	
	\begin{definition}
		An injective homomorphism $(\rho,u):\Gamma\to G\ltimes_\tR V$ is called $(Q_\tR^+,Q_\tR^-)$-\emph{affine Anosov} if it is partially $(Q_\tR^+,Q_\tR^-)$-{affine Anosov} and 
		\[0\notin\tM\text{-}\Spec(\rho,u).\]
	\end{definition}
	
	We note that the examples of proper affine action constructed by Smilga \cite{Smilga3} are examples of $(Q_\tR^+,Q_\tR^-)$-{affine Anosov} representations in $G\ltimes_\tR V$. Also, it is clear that if $(\rho,u)$ is affine Anosov with respect to $Q_\tR^\pm$, then $(\rho,u)(\Gamma)$ acts properly discontinuously on $V$. The other way round is also true but under certain``rank one" like conditions.
	
	\begin{theorem}[\cite{Ghosh7}]\label{thm.equilin}
		Suppose $\Gamma$ is word hyperbolic, $\tR$ is non-swinging and $(\rho,u):\Gamma\to G\ltimes_\tR V$ is an injective homomorphism such that $\rho$ is Anosov in $G$ with respect to $P_\tR^\pm$ and $\tM\text{-}\Spec(\rho,u)$ is a subset of some one dimensional subspace of $V^0$. Then $(\rho,u)(\Gamma)$ acts properly discontinuously on $V$ if and only if $(\rho, u)$ is affine Anosov in $G\ltimes_\tR V$ with respect to $Q_\tR^\pm$.
	\end{theorem}
	
	We note that the above theorem is a generalization of previous results by Ghosh--Treib \cite{GT} and Goldman--Labourie--Margulis \cite{GLM}.

	\section{Isospectral rigidity: Jordan and Cartan projections}
	
	We consider $\sGL_n(\C)$ and let $X$ be a proper complex algebraic subvariety of it. Suppose $X$ is the union of complex algebraic subvarieties $X_i$ for $1\leq i\leq l$ such that $X_i$'s are irreducible. They are called the \emph{irreducible components} of $X$. We consider the dimensions of $X_i$'s and denote the maximum by $\mdim(X)$. We denote the number of irreducible components of $X$ by $\irr(X)$ i.e. $\irr(X)=l$. Moreover, we denote the number of irreducible components of $X$ of the maximal dimension $\mdim(X)$ by $\irr_{\mdim}(X)$. We consider the degrees of $X_i$'s and denote the maximum among them by $\mdeg(X)$. 
	
	Suppose $x,y\in\N$. We define a function $f_{x,y}:\N\to\N\times\N$ which will be used to describe the explicit constants in the statements of the results in this section. Suppose $(x_0,y_0):= (x,y)$ and inductively we define $(x_{n+1},y_{n+1}):=(x_n^{2^{x_n}}y_n^{x_n2^{x_n}},y_n^{2^{x_n}})$. Then $f_{x,y}(n):=(x_n,y_n)$. We denote the projection onto the first coordinate of $\sum_{i=0}^{n}f_{x,y}(i)$ by $M(x,y,n)$ and define 
	\[M_{\irr,\mdeg,\mdim}(X):=M(\irr(X),\mdeg(X),\mdim(X)).\]
	
	Suppose $\Gamma$ is a finitely generated group with generating set $S$. We define:
	\[B_S(n):=\{\gamma\in\Gamma\mid |\gamma|_S\leq n\}.\]
	We also assume that $G$ is a connected semisimple real Lie group without any compact factors and with trivial center. We choose a basis of $\fg$ and obtain an isomorphism $\iota: \SL(\fg)\to \SL(n,\R)$ where $n=\dim\fg$. We note that $\iota\circ \tAd$ is a rational embedding of $G$ inside $\SL(n,\R)$ (see Proposition 4.4.5 (ii) of \cite{Sp}). We will treat $G$ as the connected component of a real algebraic subgroup of $\SL(n,\R)$. We note that $\SL(n,\C)$ is the complexification of $\SL(n,\R)$ and $\SL(n,\R)$ is complex Zariski dense inside $\SL(n,\C)$.
	We denote the complex Zariski closure inside $\SL(n,\C)$ of any algebraic subgroup $H$ of $\SL(n,\R)$ by $H_\C$. Let 
	\[\sDiag:\sGL(n,\C)\times \sGL(n,\C) \to \sGL(2n,\C)\] 
	be the standard diagonal embedding. We note that it is a regular embedding. It follows that the diagonal embedding sends $G_\C\times G_\C$ into a closed subgroup of $\sGL(2n,\C)$. 
	
	Let $\sHom(\Gamma,G)$ be the space of all injective homomorphisms of $\Gamma$ into $\sG$ and $\sHom^{\zd}(\Gamma,G)$ be the subspace of $\sHom(\Gamma,G)$ consisting of all real Zariski dense representations. Suppose $Y$ be a real algebraic subvariety of $G\times G$ such that there exist two normal subgroups $N_1,N_2$ of $G$ and a smooth isomorphism $\sigma:G/N_1\to G/N_2$ such that $Y=\{(g,h)\mid hN_2=\sigma(gN_1)\}$. We denote the collection of all such subvarieties by $\cS(G)$. Suppose $Z$ is some complex algebraic subvariety of $ \sGL(n,\C)$. We define 
	\[M_{G,Z}:=\max\{M_{\irr,\mdeg,\mdim}(\sDiag(Y_\C)\cap Z)\mid Y\in\cS(G)\}.\] We note that $M_{G,Z}$ is finite and is independent of the generating set.
	
	Suppose $g\in\sGL(2n,\C)$ and for $i,j\in\{1,2\}$ suppose $g_{i,j}$'s are square matrices of length $n$ such that 
	\[g=\begin{bmatrix} g_{11} & g_{12}\\ g_{21} & g_{22} \end{bmatrix}. \]
	
	Let $f_\kappa:\sGL(2n,\C)\to\C$ be such that $f_\kappa(g):=\ttr(g_{11}^tg_{11}-g_{22}^tg_{22})$ and $Z_{f_\kappa}$ be the zero set of $f_\kappa$.
	
	\begin{theorem}[\cite{Ghosh6}]
		Suppose $\rho_1,\rho_2\in\sHom^{\zd}(\Gamma,G)$ satisfy $\kappa(\rho_1(\gamma))=\kappa(\rho_2(\gamma))$ for all $\gamma \in B_S(M_{G,Z_{f_\kappa}})$. Then there exists an automorphism $\sigma:G \to G$ such that $\sigma \circ \rho_1=\rho_2$. 
	\end{theorem}
	
	Suppose $g$ is an element of $G$. We recall that $g$ is called \emph{loxodromic} if the Jordan projection $\tJd(g)$ lies in $\fa^{++}$. We call an element of $\sHom(\Gamma, G)$ \emph{loxodromic} if $\rho(\gamma)$ is loxodromic for all $\gamma\in\Gamma$. 
	
	Let $f_\tJd:\sGL(2n,\C)\to\C$ be such that $f_\tJd(g)=\ttr(g_{11}^2-g_{22}^2)$ and $Z_{f_\tJd}$ is the zero set of $f_\tJd$.
	
	\begin{theorem}[\cite{Ghosh6}]
		Suppose $G$ is real split and $\rho_1,\rho_2\in\sHom^{\zd}(\Gamma,G)$ are loxodromic, and  $\tJd(\rho_1(\gamma))=\tJd(\rho_2(\gamma))$
		for all $\gamma \in B_S(M_{G,Z_{f_\tJd}})$. Then there exists an automorphism $\sigma:G \to G$ such that we have $\sigma \circ \rho_1=\rho_2$.
	\end{theorem}
	
	We note that the assumption of Zariski density is crucial in both the above results. In fact, without this assumption there are counterexamples to these results. On another note, we can drop the Zariski density assumption when both the representations are in the Hitchin component (see \cite{Ghosh6}).

	\section{Isospectral rigidity: Margulis invariants}
	
	Suppose $G$ is split and $\tR:G\to\SL(V)$ is an irreducible algebraic representation. We denote the space of all representations of $\Gamma$ in $G\ltimes_\tR V$ by $\sHom(\Gamma,G\ltimes_\tR V)$ and the space of all real Zariski dense representations by $\sHom^{\zd}(\Gamma,G\ltimes_\tR V)$. We recall that $V=V^+\oplus V^0\oplus V^-$. We assume that $\dim(V^0)=q$ is non-zero and $\dim(V)=k$. We note that when $V=\fg$ and $\tR$ is the adjoint representation, the Margulis invariants are infinitesimal versions of Jordan projections (see Ghosh \cite{Ghosh7} and Sambarino \cite{Samba2}). It is natural to expect that rigidity results as that of the previous section should also hold in the affine case. Indeed, it holds and we describe the affine versions of the results from the previous section in this section.
	
	Suppose $I:V\to V$ is the identity transformation and $T:V\to V$ is any linear transformation. We recall the notion of a characteristic polynomial and to suit our purposes modify it and define the following polynomial:
	\[\chi^{red}_T(x)=\sum_{j=0}^{k-q}(-1)^{k-q-j}\ttr(\wedge^{k-q-j}(I-T))(1-x)^j.\]
	Suppose $A\in\mathfrak{gl}(k,\C)$, $v,w\in\C^k$ and $a\in\C$. We use the following notation:
	\[[A,v,w^t,a]:=\begin{bmatrix}A&v\\w^t&a\end{bmatrix}.\] 
	Moreover, for $g\in\sGL(2k+2,\C)$ we use the following notation:
	\[g=\begin{bmatrix}
		A_{11}&v_{11}&A_{12}&v_{12}\\
		w_{11}^t&a_{11}&w_{12}^t&a_{12}\\
		A_{21}&v_{21}&A_{22}&v_{22}\\
		w_{21}^t&a_{21}&w_{22}^t&a_{22}
	\end{bmatrix}.\]

	Let $f_\tM:\sGL(k+1,\C)\to\C^k$ be such that $f_\tM([A,v,w^t,a]):=\chi^{red}_A(A)v$ and $Z_{f_\tM}$ be the zero set of $f_\tM$.
	\begin{theorem}[\cite{Ghosh6}]
		Suppose $(\rho,u),(\varrho,v)\in\sHom(\Gamma,G\ltimes_\tR V)$ are such that $\rho$ is Zariski dense and loxodromic and $\varrho$ is a conjugate of $\rho$. Also, suppose $\tM_{\rho,u}(\gamma)=\tM_{\varrho,v}(\gamma)$ for all $\gamma \in B_S(M_{\irr,\mdeg,\mdim}((G_\C\ltimes_\tR V_\C)\cap Z_{f_\tM}))$. Then $(\rho,u)$ and $(\varrho,v)$ are also conjugate to each other.
	\end{theorem}
	Henceforth, we assume that $\tR$ preserves a non-degenerate symmetric real bilinear form on $V$. We note that when $V=\fg$ and $\tR$ is the adjoint representation, the Killing form is a non-degenerate symmetric real bilinear form on $\fg$ preserved by $\tAd$. Let $\tB_\tR$ be the complex bilinear form on $V_\C$ corresponding to the real bilinear form on $V$. We abuse notation and for simplicity denote $\tB_\tR(v,v)$ by $\tB_\tR(v)$. 
	
	Let $F:\sGL(2k+2,\C)\to\C$ be such that 
	\[F(g):=\tB_\tR(\chi^{red}_{A_{11}}(1)\chi^{red}_{A_{22}}(A_{22})v_{22})-\tB_\tR(\chi^{red}_{A_{22}}(1)\chi^{red}_{A_{11}}(A_{11})v_{11})\] 
	and $f:\sGL(k+1,\C)\to\C$ be such that $f([A,v,w^t,a]):=\tB_\tR(\chi^{red}_A(A)v)$. We denote the zero set of $F$ and $f$ respectively by $Z_F$ and $Z_f$. We denote the maximum of $M_{G\ltimes_\tR V,Z_F}$ and $M_{\irr,\mdeg,\mdim}((G_\C\ltimes_\tR V_\C)\cap Z_f)$ by $m$.
	
	\begin{theorem}[\cite{Ghosh6}]
		Suppose  $(\rho,u),(\varrho,v)\in\sHom(\Gamma,G\ltimes_\tR V)$ are such that both $\rho$ and $\varrho$ are Zariski dense and loxodromic representations. Moreover, suppose $\tB_\tR(\tM_{\rho,u}(\gamma))=\tB_\tR(\tM_{\varrho,v}(\gamma))$ for all $\gamma \in B_S(m)$. Then the following holds:
		\begin{enumerate}
			\item If $\tB_\tR(\tM_{\rho,u}(\gamma))=0$ for all $\gamma\in B_S(m)$, then the real Zariski closure of both $(\rho,u)(\Gamma)$ and $(\varrho,v)(\Gamma)$ are conjugate to $G$.
			\item If $\tB_\tR(\tM_{\rho,u}(\eta))\neq0$ for some $\eta\in B_S(m)$, then $(\rho,u)$ and $(\varrho,v)$ are conjugate via some element of $\sGL(V)\ltimes V$.
		\end{enumerate}
	\end{theorem}
	
	We note that the above results imply that two representations which gives rise to proper affine actions and whose Margulis invariants match up to some well chosen finite set are same up to automorphisms of $G\ltimes_\tR V$. In fact, the same holds for affine Anosov representations too. These results substantially generalize previous results by Drumm--Goldman \cite{DG}, Charette--Drumm \cite{CD}, Kim \cite{Kim} and Ghosh \cite{Ghosh5}.
	
	\section{The case of split pseudo-orthogonal groups}
	
	The theory of affine deformations of Anosov representations is more well developed in the case of the split orthogonal groups $\mathsf{SO}(n,n-1)$ acting on $\R^{2n-1}$. In this section we present known results for this particular group. We expect that many of these results will continue to hold for more general settings. The generalizations are a work in progress.
	
	We recall that there exists non-abelian free subgroups of $\mathsf{SO}(n,n-1)\ltimes\R^{2n-1}$ acting properly discontinuously on $\R^{2n-1}$ only when $n$ is even. We standardize our notations. We denote the identity matrix of length $2n$ by $I_{2n}$ and consider the following quadratic form on $\R^{4n}$: 
	\[I_{2n,2n}\defeq
	\begin{bmatrix}
		I_{2n} & 0\\
		0 & -I_{2n}
	\end{bmatrix}.\]
	We denote the group of all those linear transformations of $\R^{4n}$ which preserve $I_{2n,2n}$ by $\mathsf{SO}(2n,2n)$ and its connected component containing the identity by $\mathsf{SO}_0(2n,2n)$. We denote the standard basis of $\R^{4n}$ by $e_1,\dots,e_{4n}$ and note that the group formed by all those elements of $\mathsf{SO}(2n,2n)$ which fixes $e_{4n}$ is a copy of $\mathsf{SO}_0(2n,2n-1)$. We denote the corresponding Lie algebras by $\mathfrak{so}(2n,2n)$ and $\mathfrak{so}(2n,2n-1)$. We note that $\mathfrak{so}(2n,2n-1)e_{4n}=\{0\}$ and $\mathfrak{so}(2n,2n)e_{4n}=\R^{4n-1}\times\{0\}\subset\R^{4n}$. We call a subspace of $\R^{4n-1}$ a \emph{null subspace} if it is orthogonal to some maximal isotropic subspace of $\R^{4n-1}$. The following theorem shows that affine Anosov representations in $\mathsf{SO}(2n,2n-1)\ltimes\R^{4n-1}$ appear due to the following special phenomenon: all elements of $\mathsf{SO}(2n,2n)$ which are loxodromic in $\SL(4n,\R)$ are also loxodromic in $\mathsf{SO}(2n,2n)$ but all loxodromic elements of $\mathsf{SO}(2n,2n)$ need not be loxodromic in $\SL(4n,\R)$. 
	
	\begin{theorem}[\cite{Ghosh3}]
		Let $\{\rho_t\}_{t\in(-1,1)}$ be an analytic one parameter family of representations of $\Gamma$ in $\mathsf{SO}(2n,2n)$ with $\rho_0(\Gamma)\subset\mathsf{SO}(2n,2n-1)$. Let ${U}$ be the tangent vector to $\{\rho_t\}_{t\in(-1,1)}$ at $\rho=\rho_0$ and ${u}={U}e_{4n}$. Suppose $(\rho,u)$ is affine Anosov with respect to the stabilizer of a null subspace of $\R^{4n-1}$. Then there exists $\epsilon>0$ such that for all $t$ with $|t|\in(0,\epsilon)$, $\rho_t$ is Anosov in $\mathsf{SL}(4n,\R)$ with respect to the stabilizer of an $n$-dimensional subspace.
	\end{theorem}
	
	We denote the space of conjugacy classes of representations of $\Gamma$ inside $\mathsf{SO}(2n,2n-1)\ltimes\R^{4n-1}$ which are affine Anosov with respect to the stabilizer of some null subspace of $\R^{4n-1}$ by $\cA_n$. Henceforth, we will replace the phrase ``affine Anosov with respect to the stabilizer of some null subspace of $\R^{4n-1}$" by the expression ``affine Anosov" for simplicity. We denote the conjugacy classes of representations of $\Gamma$ inside $\mathsf{SO}(2n,2n-1)$ which are linear parts of such affine Anosov representations by $\cL_n$. Moreover, let $\cP_n$ be the collection of conjugacy classes of all affine representations in $\sHom(\Gamma,\mathsf{SO}(2n,2n-1)\ltimes\R^{4n-1})$ whose linear parts are the linear parts of some affine Anosov representations. We note that $\cP_n$ is the space of conjugacy classes of all partially affine Anosov representations. Clearly, we have $\cA_n\subset\cP_n$. We note that $\cL_n$, $\cP_n$ and $\cA_n$ are all open subsets of the respective representation spaces. Moreover, we know that $\cA_n$ is fibered over $\cL_n$ with fibers which are a disjoint union of two convex cones which differ by a sign. 
	
	We note that in the case of split pseudo-orthogonal groups $\mathsf{SO}(2n,2n-1)$, the Margulis invariants lie in a line. Hence, Theorem \ref{thm.equilin} is applicable in this context. Moreover, in these cases Margulis invariants can be treated as real numbers. We denote the set of periodic orbits of the Gromov flow space by $O$. Suppose $(\rho,u)$ is an affine Anosov representation. We define, $R_T(\rho,u)\defeq\{\gamma\in O\mid\tM_{\rho,u}(\gamma)\leqslant T\}$ and note that the following limit exists and is positive (see Ghosh \cite{Ghosh4}):
	\[h_{\rho,u}=\lim_{T\to\infty}\frac{1}{T}\log|R_T(\rho,u)|.\]
	The quantity $h_{\rho,u}$ is called the \emph{topological entropy} of $(\rho,u)$. The topological entropy satisfies the following properties for $c>0$:
	\begin{enumerate}
		\item $h_{\rho,cu}=\frac{1}{c}h_{\rho,u}$,
		\item if $h_{\rho,u}=c=h_{\rho,v}$, then $h_{\rho,(1-t)u+tv}< c$ for all $0\leq t\leq1$. 
	\end{enumerate}
	
	Also, for a partially affine Anosov representation $(\varrho,v)$ we note that the following limit exists (see Ghosh \cite{Ghosh4}):
	\[I((\rho,u),(\varrho,v))=\lim_{T\to\infty}\frac{1}{|R_T(\rho,u)|}\sum_{\gamma\in R_T(\rho,u)}\frac{\tM_{\varrho,v}(\gamma)}{\tM_{\rho,u}(\gamma)}.\]
	The quantity $I$ is called the \emph{intersection number} of $(\rho,u)$ and $(\varrho,v)$. As the Margulis invariants are preserved under conjugacy, we obtain well defined maps $h:\cA_n\to\R$ and $I:\cA_n\times\cP_n\to\R$. We note that the maps $h,I$ are real analytic (see Ghosh \cite{Ghosh4}).
	
	Suppose $\{(\rho_t,u_t)\}_{t\in\R}$ is an analytic one parameter family of affine Anosov representations such that $(\rho_0,u_0)=(\rho,u)$ and for all $\gamma\in\Gamma$:
	$${U}(\gamma)=\left.\frac{d}{dt}\right|_{t=0}(\rho_t,u_t)(\gamma)(\rho,u)(\gamma)^{-1}.$$ 
	We note that $[U]\in\mathsf{H}^1_{\sad\circ(\rho,u)}(\Gamma,\mathsf{so}(2n,2n-1)\oplus\R^{4n-1})\cong \sT_{[\rho,u]}\cA_n$. We define the \emph{pressure form} $\tP:\sT\cA_n\times\sT\cA_n\to\R$ as follows:
	\[\tP_{[\rho,u]}([U],[U])\defeq\left.\frac{d^2}{dt^2}\right|_{t=0}\frac{h_{\rho_t}}{h_\rho}I(\rho,\rho_t).\]	
	We note that $\tP$ is well defined and is a positive semi-definite bilinear form. Suppose $\cA_{n,k}$ are the constant entropy sections of $\cA_n$ with entropy $k$ i.e.
	\[\cA_{n,k}:=\{\rho\in\cA_n\mid h_\rho=k\}.\]
	We note that $\cA_{n,k}$ is an analytic submanifold of $\cA_n$ of codimension one. 
	\begin{theorem}[\cite{Ghosh4}]
		The restriction of the pressure form to the constant entropy sections of $\cA_n$ i.e. $\tP:\sT\cA_{n,k}\times\sT\cA_{n,k}\rightarrow\R$ is a Riemannian metric for all $k>0$ and the pressure form on $\cA_n$ is positive semi-definite with rank $\dim(\cA_n)-1$. 
	\end{theorem}

	\bibliography{Library.bib}
	\bibliographystyle{alpha}
	
\end{document}